\documentclass[12pt]{article}
\usepackage{amsfonts}

\usepackage{graphicx}
\usepackage{amsmath}

\newtheorem{theorem}{Theorem}

\newtheorem{conjecture}[theorem]{Conjecture}
\newtheorem{corollary}[theorem]{Corollary}

\newtheorem{lemma}[theorem]{Lemma}

\newtheorem{proposition}[theorem]{Proposition}
\newtheorem{remark}[theorem]{Remark}

\newenvironment{proof}[1][Proof]{\noindent\textbf{#1.} }{\hfill $\Box $}

\begin{document}

\title{Logarithmic vector fields and hyperbolicity}
\author{Erwan Rousseau\thanks{\textit{Mathematics Subject Classification (2000)}:
Primary: 32Q45, 14J70; \textit{Key words}: Kobayashi hyperbolicity; jet differentials; logarithmic vector fields.}}
\date{}
\maketitle

\begin{abstract}
Using vector fields on logarithmic jet spaces we obtain some new positive results for the logarithmic Kobayashi conjecture about the hyperbolicity of complements of curves in the complex projective plane. We are interested here in the cases where logarithmic irregularity is strictly smaller than the dimension. In this setting, we study the case of a very generic curve with two components of
degrees $d_{1}\leq d_{2}$ and prove the hyperbolicity of the
complement if the degrees satisfy either $d_{1}\geq4,$ or $d_{1}=3$ and $d_{2}\geq 5$, or $d_{1}=2$ and $d_{2}\geq
8$, or $d_{1}=1$ and $d_{2}\geq 11$. We also prove that the complement of a very generic
curve of degree $d$ at least equal to 14 in the complex projective
plane is hyperbolic, improving slightly, with a new proof, the former bound obtained by El Goul.
\end{abstract}

\section{Introduction}

A complex manifold $X$ is hyperbolic in the sense of S. Kobayashi if the
hyperbolic pseudodistance defined on $X$ is a distance (see, for example,
\cite{Ko98}). In the case of hypersurfaces in $\mathbb{P}^{n}$, we have the
logarithmic Kobayashi conjecture \cite{Ko70}:

\begin{conjecture}
\label{c1}$\mathbb{P}^{n}\backslash X$ $(n\geq 2)$ is hyperbolic for a
generic hypersurface $X\subset \mathbb{P}^{n}$ of degree $\deg X\geq 2n+1.$
\end{conjecture}

Here we will study the case of complements of curves in $\mathbb{P}^{2}.$
Several authors have studied this case, especially when the curve has several
irreducible components. It is well known that the conjecture is the more difficult the smaller the logarithmic irregularity (equivalently the number of irreducible components) is. The conjecture is known to be true for logarithmic irregularity equal to 2 or more (equivalently 3 or more irreducible components). We refer to \cite{DSW1}, \cite{DSW2} and \cite{BeDu} for the details and the references for these cases. In the case of logarithmic irregularity equal to 2, Dethloff and Lu proved in \cite{DL04} that every Brody curve in the complement of a normal crossing curve in $\mathbb{P}^{2}$ of degree at least $4$ consisting of three components is algebraically degenerate. See also \cite{NWY} for more general results on the algebraic degeneracy of entire curves when the logarithmic irregularity is equal to the dimension of the manifold.

When the logarithmic irregularity is strictly smaller than the dimension of the manifold, much less is known. In this paper, we are interested in the more difficult cases where the logarithmic irregularity is strictly smaller than $2$, i.e the curve is either smooth or has two irreducible components. For the complement of smooth curves, studying
the compact analogue of the above conjecture, Demailly and El Goul obtained
in \cite{DEG00} that complements of very generic curves in $\mathbb{P}^{2}$ of degree $d\geq 21$ are hyperbolic. Later, using logarithmic jets, El Goul improved that result in \cite{E.G} obtaining the bound 15. Using different techniques we obtain here

\begin{theorem}
\label{t1}Let $C$ be a very generic irreducible complex algebraic curve in $\mathbb{P}%
^{2}$ of degree $d$. Then $\mathbb{P}^{2}\backslash C$ is hyperbolic and
hyperbolically embedded in $\mathbb{P}%
^{2}$ if $d\geq14$.
\end{theorem}

Previously in \cite{Rou03} we obtained some results for the two-components case. Using the same techniques as in the proof of the previous theorem we improve them by the following result

\begin{theorem}
\label{t2}Let $C=C_{1}\cup C_{2}$ be a very generic complex algebraic curve in $\mathbb{P}%
^{2}$ having two irreducible components $C_{1}$ and $C_{2}$ of degrees $d_{1}\leq d_{2}$. Then $\mathbb{P}^{2}\backslash C$ is hyperbolic and
hyperbolically embedded in $\mathbb{P}%
^{2}$ if the degrees satisfy
$$\begin{array}{cccc}
either & d_{1} \geq 4, \\
or & d_{1} =3 & and & d_{2}\geq 5,\\
or & d_{1}=2 &and& d_{2}\geq 8,\\
or & d_{1}=1 &and& d_{2}\geq 11.
\end{array}
$$
\end{theorem}

The proofs of these two results are based on techniques  introduced by Siu and Paun (see \cite{SY04}, \cite{Paun05}, \cite{Rou06} and \cite{Rou06bis}).

The first one is a generalization in the logarithmic setting of an approach
initiated by Clemens \cite{Cle}, Ein \cite{Ein}, Voisin \cite{Voi} and used
by Y.-T. Siu \cite{SY04} and M. Paun \cite{Paun05} to construct vector fields on the total space of
hypersurfaces in the projective space. Here we construct vector fields on
logarithmic spaces.

The second one is based on bundles of logarithmic jet differentials (see
\cite{DL96}). The idea, in hyperbolicity questions, is that global sections
of these bundles vanishing on ample divisors provide algebraic differential
equations for any entire curve $f:\mathbb{C}\rightarrow X\backslash D$ where
$D$ is a normal crossing divisor on $X.$ Therefore, the main point is to
produce enough algebraically independent global holomorphic logarithmic jet
differentials.

\section{Logarithmic jet bundles}

In this section we recall briefly the basic facts and results of J. Noguchi in \cite
{No86} about logarithmic jet bundles. We refer to \cite{DL96} and \cite{Rou06bis} for details.

Let $X$ be a complex manifold of dimension $n$. Denote by $J_{k}X$ the $k$-jet bundle over $X$. Let $T_{X}^{\ast }$ be the holomorphic cotangent bundle over $X.$ Take a
holomorphic section $\omega \in H^{0}(O,T_{X}^{\ast })$ for some open subset
$O.$ For $j_{k}(f)\in J_{k}X_{\left| O\right. },$ we have $f^{\ast }\omega
=Z(t)dt$ and a well defined holomorphic mapping
\begin{equation*}
\widetilde{\omega }:J_{k}X_{\left| O\right. }\rightarrow \mathbb{C}%
^{k};j_{k}(f)\rightarrow \left( \frac{d^{j}Z}{dt^{j}}(0)\right) _{0\leq
j\leq k-1}.
\end{equation*}

If $\omega _{1},\dots,\omega _{n}$ are holomorphic 1-forms on $O$
such that $\omega _{1} \wedge \dots \wedge \omega _{n}$ vanishes
nowhere, then we have a biholomorphic map
\begin{equation*}
\left( \widetilde{\omega _{1}},\dots,\widetilde{\omega _{n}}\right) \times \pi
:J_{k}X_{\left| O\right. }\rightarrow \left( \mathbb{C}^{k}\right)
^{n}\times O,
\end{equation*}
which gives the trivialization of $J_{k}X_{\left| O\right.}$ associated to $\omega _{1},\dots,\omega _{n}.$

Let $\overline{X}$ be a complex manifold with a normal crossing divisor $D.$
Consider the log manifold $(\overline{X},D)$. Let $X=\overline{X}%
\backslash D.$ Denote by $\overline{T}_{X}^{\ast }=T_{\overline{X}%
}^{\ast }(\log D)$ the logarithmic cotangent sheaf. 

Let $s\in H^{0}(O,J_{k}\overline{X})$ be a holomorphic section over an open
subset $O\subset \overline{X}.$ We say that $s$ is a logarithmic $k$-jet
field if the map $\widetilde{\omega }\circ s_{\left| O^{\prime }\right.
}:O^{\prime }\rightarrow \mathbb{C}^{k}$ is holomorphic for all $\omega \in
H^{0}(O^{\prime },\overline{T}_{X}^{\ast })$ for all open subsets $O^{\prime
}$ of $O.$ The set of logarithmic $k$-jet fields over open subsets of $%
\overline{X}$ defines a subsheaf of the sheaf $J_{k}\overline{X},$ which we
denote by $\overline{J}_{k}X.$ $\overline{J}_{k}X$ is the sheaf of sections
of a holomorphic fibre bundle over $\overline{X},$ denoted again $\overline{J%
}_{k}X$ and called the logarithmic $k-$jet bundle of $(\overline{X},D).$

A log-morphism $F:(\overline{X}^{\prime },D^{\prime })\rightarrow (\overline{%
X},D)$ induces a canonical map
\begin{equation*}
F_{k}:\overline{J}_{k}X^{\prime }\rightarrow \overline{J}_{k}X.
\end{equation*}

We can express the local triviality of $\overline{J}_{k}X$ explicitly in
terms of coordinates. Let $z_{1,}\dots,z_{n}$ be coordinates in an open set $%
O\subset \overline{X}$ in which $D=\{z_{1}z_{2}\dots z_{l}=0\}.$ Let $\omega
_{1}=\frac{dz_{1}}{z_{1}},\dots\omega _{l}=\frac{dz_{l}}{z_{l}},\omega
_{l+1}=dz_{l+1},\dots,\omega _{n}=dz_{n}.$ Then we have a biholomorphic map
\begin{equation*}
\left( \widetilde{\omega _{1}},\dots,\widetilde{\omega _{n}}\right) \times \pi
:\overline{J}_{k}X_{\left| O\right. }\rightarrow \left( \mathbb{C}%
^{k}\right) ^{n}\times O.
\end{equation*}

Let $s\in H^{0}(O,\overline{J}_{k}X)$ be given by $s(x)=(\xi
_{j}^{(i)}(x),x) $ in this trivialization where the indices $i$ correspond
to the orders of derivative. Then the same $s$ considered as an element of $%
H^{0}(O,J_{k}\overline{X})$ and trivialized by $\omega
_{1}=dz_{1},\dots\omega _{n}=dz_{n}$ is given by $s(x)=(\widehat{\xi }%
_{j}^{(i)}(x),x)$ where
\begin{equation*}
\widehat{\xi }_{j}^{(i)}=\left\{
\begin{array}{c}
z_{i}(\xi _{j}^{(i)}+g_{i}(\xi _{j}^{(1)},\dots,\xi _{j}^{(i-1)})):j\leq l, \\
\xi _{j}^{(i)}:j\geq l+1.
\end{array}
\right.
\end{equation*}

The $g_{i}$ are polynomials in the variables $\xi _{j}^{(1)},\dots,\xi
_{j}^{(i-1)},$ obtained by expressing first the different components $\xi
_{j}^{(i)}$ of $\left( \widetilde{\frac{dz_{i}}{z_{i}}}\right) \circ s(x)$
in terms of the components $\widehat{\xi }_{j}^{(i)}$ of \ the components $%
\widehat{\xi }_{j}^{(i)}$ of $\widetilde{dz_{i}}\circ s(x)$ by using the
chain rule, and then by inverting this system.

In summary, we have a holomorphic coordinate system on $\overline{J}_{k}X_{\left| O\right. }$ given by $(\xi_{1}^{(1)}\dots,\xi_{n}^{(k)};z_{1},\dots,z_{n})$ and one on  $J_{k}\overline{X}_{\left| O\right. }$ given by $(\widehat{\xi }_{1}^{(1)}\dots,\widehat{\xi }_{n}^{(k)};z_{1},\dots,z_{n})$. The previous relation exhibits the sheaf inclusion $\overline{J}_{k}X_{\left| O\right. } \subset J_{k}\overline{X}_{\left| O\right. }$.
We will use these coordinates for the computations of the next section.

\section{Logarithmic vector fields}

\subsection{The smooth case}
In this section we generalize the approach used in \cite{Paun05} (see also 
\cite{PaRou} and \cite{Rou06bis}) to logarithmic jet bundles. Once we have the logarithmic tools of the previous section, the explicit construction of the vector fields on logarithmic jet spaces is very similar to the compact case treated in \cite{Paun05}. Nevertheless, we provide the details of the computations for the convenience of the reader.

Let $\mathcal{X}\subset \mathbb{P}^{2}\times \mathbb{P}^{N_{d}}$ be the
universal curve of degree $d$ given by the equation
\begin{equation*}
\underset{\left| \alpha \right| =d}{\sum }a_{\alpha }Z^{\alpha }=0,\text{
where }[a]\in \mathbb{P}^{N_{d}}\text{ and }[Z]\in \mathbb{P}^{2}.
\end{equation*}

We use the
notations: for $\alpha =(\alpha _{0},\dots,\alpha _{2})\in \mathbb{N}^{3},$ $%
\left| \alpha \right| =\sum_{i}\alpha _{i}$ and if $Z=(Z_{0},Z_{1},Z_{2})$
are homogeneous coordinates on $\mathbb{P}^{2},$ then $Z^{\alpha }=\prod
Z_{j}^{\alpha _{j}}.$ $\mathcal{X}$ is a smooth hypersurface of degree $%
(d,1) $ in $\mathbb{P}^{2}\times \mathbb{P}^{N_{d}}.$

We consider the log-manifold $(\mathbb{P}^{2}\times \mathbb{P}^{N_{d}},%
\mathcal{X)}$. We denote by $\overline{J_{2}}(\mathbb{P}^{2}\times \mathbb{P}%
^{N_{d}})$ the manifold of the logarithmic 2-jets$,$ and $\overline{J_{2}^{v}%
}(\mathbb{P}^{2}\times \mathbb{P}^{N_{d}})$ the submanifold of $\overline{%
J_{2}}(\mathbb{P}^{2}\times \mathbb{P}^{N_{d}})$ consisting of 2-jets
tangent to the fibers of the projection $\mathbb{P}^{2}\times \mathbb{P}%
^{N_{d}}\rightarrow \mathbb{P}^{N_{d}}.$

We are going to construct meromorphic vector fields on $\overline{J_{2}^{v}}(%
\mathbb{P}^{2}\times \mathbb{P}^{N_{d}}).$

Let us consider
\begin{equation*}
\mathcal{Y}=(a_{000d}Z_{3}^{d}+\underset{\left| \alpha \right| =d}{\sum }%
a_{\alpha }Z^{\alpha }=0)\subset \mathbb{P}^{3}\times U,
\end{equation*}
where $U:=(a_{000d}\neq 0)\cap \left( \underset{\left| \alpha \right|
=d,\alpha _{n+2}=0}{\cup }(a_{\alpha }\neq 0)\right) \subset \mathbb{P}%
^{N_{d}+1}.$ We have the projection $\pi :\mathcal{Y}\rightarrow \mathbb{P}%
^{2}\times \mathbb{P}^{N_{d}}$ and $\pi ^{-1}(\mathcal{X})=(Z_{3}=0)$ $:=H.$
Therefore we obtain a log-morphism $\pi :(\mathcal{Y},H)\rightarrow (\mathbb{%
P}^{2}\times \mathbb{P}^{N_{d}},\mathcal{X})$ which induces a dominant map
\begin{equation*}
\pi _{2}:\overline{J_{2}^{v}}(\mathcal{Y})\rightarrow \overline{J_{2}^{v}}(%
\mathbb{P}^{2}\times \mathbb{P}^{N_{d}}).
\end{equation*}

Let us consider the set $\Omega _{0}:=(Z_{0}\neq 0)\times (a_{000d}\neq
0)\subset \mathbb{P}^{3}\times U.$ We assume that global coordinates are
given on $\mathbb{C}^{3}$ and $\mathbb{C}^{N_{d}+1}.$ The equation of $%
\mathcal{Y}$ becomes
\begin{equation*}
\mathcal{Y}_{0}:=(z_{3}^{d}+\underset{\alpha }{\sum }a_{\alpha }z^{\alpha
}=0).
\end{equation*}

Following \cite{DL96} as explained above, we can obtain explicitly a
trivialization of $\overline{J_{2}}(\Omega _{0}).$ Let $\omega
^{1}=dz_{1},\omega ^{2}=dz_{2},\omega ^{3}=\frac{dz_{3}}{z_{3}}.$ Then we
have a biholomorphic map
\begin{equation*}
\overline{J_{2}}(\Omega _{0})\rightarrow \mathbb{C}^{3}\times U\times
\mathbb{C}^{3}\times \mathbb{C}^{3},
\end{equation*}
where the coordinates will be denoted $(z_{i},a_{\alpha },\xi _{j}^{(i)}).$

Let us write the equations of $\overline{J_{2}^{v}}(\mathcal{Y}_{0})$ in this
trivialization. We have $\overline{J_{2}^{v}}(\mathcal{Y}_{0})=J_{2}^{v}(%
\mathcal{Y}_{0})\cap \overline{J_{2}}(\Omega _{0}).$ The equations of $%
J_{2}^{v}(\mathcal{Y}_{0})$ in the trivialization of $J_{2}(\Omega _{0})$
given by $\widehat{\omega }^{1}=dz_{1},\widehat{\omega }^{2}=dz_{2},\widehat{%
\omega }^{3}=dz_{3}\mathbb{\ }$can be written in $\mathbb{C}^{3}\times
U\times \mathbb{C}^{3}\times \mathbb{C}^{3}$ with coordinates $%
(z_{i},a_{\alpha },\widehat{\xi }_{j}^{(i)})$ as follows:
\begin{equation*}
z_{3}^{d}+\underset{\left| \alpha \right| \leq d}{\sum }a_{\alpha }z^{\alpha
}=0,
\end{equation*}
\begin{equation*}
dz_{3}^{d-1}\widehat{\xi }_{3}^{(1)}+\underset{j=1}{\overset{2}{\sum }}%
\underset{\left| \alpha \right| \leq d}{\sum }a_{\alpha }\frac{\partial
z^{\alpha }}{\partial z_{j}}\widehat{\xi }_{j}^{(1)}=0,
\end{equation*}
\begin{equation*}
dz_{3}^{d-1}\widehat{\xi }_{3}^{(2)}+d(d-1)z_{3}^{d-2}\left( \widehat{\xi }%
_{3}^{(1)}\right) ^{2}+\underset{j=1}{\overset{2}{\sum }}\underset{\left|
\alpha \right| \leq d}{\sum }a_{\alpha }\frac{\partial z^{\alpha }}{\partial
z_{j}}\widehat{\xi }_{j}^{(2)}+\underset{j,k=1}{\overset{2}{\sum }}\underset{%
\left| \alpha \right| \leq d}{\sum }a_{\alpha }\frac{\partial ^{2}z^{\alpha }%
}{\partial z_{j}\partial z_{k}}\widehat{\xi }_{j}^{(1)}\widehat{\xi }%
_{k}^{(1)}=0.
\end{equation*}

The relations between the two systems of coordinates can be computed as
explained above and are given by
\begin{eqnarray*}
\widehat{\xi }_{j}^{(i)} &=&\xi _{j}^{(i)}\text{ for }j\leq 2, \\
\widehat{\xi }_{3}^{(1)} &=&z_{3}\xi _{3}^{(1)}, \\
\widehat{\xi }_{3}^{(2)} &=&z_{3}\left(\xi _{3}^{(2)}+\left( \xi
_{3}^{(1)}\right) ^{2}\right).
\end{eqnarray*}

Therefore, to obtain the equations of $\overline{J_{2}^{v}}(\mathcal{Y}_{0})$
in the first trivialization, we just have to substitute the previous
relations
\begin{equation}
z_{3}^{d}+\underset{\left| \alpha \right| \leq d}{\sum }a_{\alpha }z^{\alpha
}=0,
\end{equation}

\begin{equation}
dz_{3}^{d}\xi _{3}^{(1)}+\underset{j=1}{\overset{2}{\sum }}\underset{\left|
\alpha \right| \leq d}{\sum }a_{\alpha }\frac{\partial z^{\alpha }}{\partial
z_{j}}\xi _{j}^{(1)}=0,
\end{equation}

\begin{equation}
dz_{3}^{d}\xi _{3}^{(2)}+d^{2}z_{3}^{d}\left( \xi _{3}^{(1)}\right) ^{2}+%
\underset{j=1}{\overset{2}{\sum }}\underset{\left| \alpha \right| \leq d}{%
\sum }a_{\alpha }\frac{\partial z^{\alpha }}{\partial z_{j}}\xi _{j}^{(2)}+%
\underset{j,k=1}{\overset{2}{\sum }}\underset{\left| \alpha \right| \leq d}{%
\sum }a_{\alpha }\frac{\partial ^{2}z^{\alpha }}{\partial z_{j}\partial z_{k}%
}\xi _{j}^{(1)}\xi _{k}^{(1)}=0.
\end{equation}

Following the method used in \cite{Paun05} and \cite{Rou06} for the compact case, we are going
to prove that $T_{\overline{J_{2}^{v}}(\mathcal{Y})}\otimes \mathcal{O}_{%
\mathbb{P}^{4}}(c)\otimes \mathcal{O}_{\mathbb{P}^{N_{d}+1}}(\ast )$ is
generated by its global sections on $\overline{J_{2}^{v}}(\mathcal{Y}%
)\backslash (\Sigma \cup p^{-1}(H)),$ where $p:$ $\overline{J_{2}^{v}}(%
\mathcal{Y})\rightarrow $ $\mathcal{Y}$ is the natural projection, $\Sigma $
a subvariety that will be defined below, and $c\in \mathbb{N}$ a constant
independant of $d.$ Consider a vector field
\begin{equation*}
V=\underset{\left| \alpha \right| \leq d}{\sum }v_{\alpha }\frac{\partial }{%
\partial a_{\alpha }}+\underset{j}{\sum }v_{j}\frac{\partial }{\partial z_{j}%
}+\underset{j,k}{\sum }w_{j}^{(k)}\frac{\partial }{\partial \xi _{j}^{(k)}},
\end{equation*}
on $\mathbb{C}^{3}\times U\times \mathbb{C}^{3}\times \mathbb{C}^{3}.$ The
conditions to be satisfied by $V$ to be tangent to $\overline{J_{2}^{v}}(%
\mathcal{Y}_{0})$ are the following
\begin{equation}
\underset{\left| \alpha \right| \leq d}{\sum }v_{\alpha }z^{\alpha }+%
\underset{j=1}{\overset{2}{\sum }}\underset{\left| \alpha \right| \leq d}{%
\sum }a_{\alpha }\frac{\partial z^{\alpha }}{\partial z_{j}}%
v_{j}+dz_{3}^{d-1}v_{3}=0,
\end{equation}

\begin{eqnarray}
\underset{j=1}{\overset{2}{\sum }}\underset{\left| \alpha \right| \leq d}{%
\sum }v_{\alpha }\frac{\partial z^{\alpha }}{\partial z_{j}}\xi _{j}^{(1)}+%
\underset{j,k=1}{\overset{2}{\sum }}\underset{\left| \alpha \right| \leq d}{%
\sum }a_{\alpha }\frac{\partial ^{2}z^{\alpha }}{\partial z_{j}\partial z_{k}%
}v_{j}\xi _{k}^{(1)}+\underset{j=1}{\overset{2}{\sum }}\underset{\left|
\alpha \right| \leq d}{\sum }a_{\alpha }\frac{\partial z^{\alpha }}{\partial
z_{j}}w_{j}^{(1)} &&  \notag \\
+d^{2}z_{3}^{d-1}v_{3}\xi _{3}^{(1)}+dz_{3}^{d}w_{3}^{(1)}=0, &&
\end{eqnarray}
\begin{eqnarray}
\underset{\left| \alpha \right| \leq d}{\sum }(\underset{j=1}{\overset{2}{%
\sum }}\frac{\partial z^{\alpha }}{\partial z_{j}}\xi _{j}^{(2)}+\underset{%
j,k=1}{\overset{2}{\sum }}\frac{\partial ^{2}z^{\alpha }}{\partial
z_{j}\partial z_{k}}\xi _{j}^{(1)}\xi _{k}^{(1)})v_{\alpha } &&  \notag \\
+\underset{j=1}{\overset{2}{\sum }}\underset{\left| \alpha \right| \leq d}{%
\sum }a_{\alpha }(\underset{k=1}{\overset{2}{\sum }}\frac{\partial
^{2}z^{\alpha }}{\partial z_{j}\partial z_{k}}\xi _{k}^{(2)}+\underset{k,l=1%
}{\overset{2}{\sum }}\frac{\partial ^{3}z^{\alpha }}{\partial z_{j}\partial
z_{k}\partial z_{l}}\xi _{k}^{(1)}\xi _{l}^{(1)})v_{j} &&  \notag \\
+\underset{\left| \alpha \right| \leq d}{\sum }(\underset{j,k=1}{\overset{2}{%
\sum }}a_{\alpha }\frac{\partial ^{2}z^{\alpha }}{\partial z_{j}\partial
z_{k}}(w_{j}^{(1)}\xi _{k}^{(1)}+w_{k}^{(1)}\xi _{j}^{(1)})+\underset{j=1}{%
\overset{2}{\sum }}a_{\alpha }\frac{\partial z^{\alpha }}{\partial z_{j}}%
w_{j}^{(2)}) &&  \notag \\
+v_{3}d^{2}z_{3}^{d-1}(\xi _{3}^{(2)}+d\left( \xi _{3}^{(1)}\right)
^{2})+2d^{2}z_{3}^{d}w_{3}^{(1)}\xi _{3}^{(1)}+dz_{3}^{d}w_{3}^{(2)}=0. &&
\end{eqnarray}

We can introduce the first package of vector fields tangent to $\overline{%
J_{2}^{v}}(\mathcal{Y}_{0}).$ We denote by $\delta _{j}\in \mathbb{N}^{2}$
the multi-index whose j-component is equal to 1 and the other are zero.

For $\alpha _{1}\geq 3:$%
\begin{equation*}
V_{\alpha }^{30}:=\frac{\partial }{\partial a_{\alpha }}-3z_{1}\frac{%
\partial }{\partial a_{\alpha -\delta _{1}}}+3z_{1}^{2}\frac{\partial }{%
\partial a_{\alpha -2\delta _{1}}}-z_{1}^{3}\frac{\partial }{\partial
a_{\alpha -3\delta _{1}}}.
\end{equation*}

For $\alpha _{1}\geq 2,\alpha _{2}\geq 1:$%
\begin{eqnarray*}
V_{\alpha }^{21} &:&=\frac{\partial }{\partial a_{\alpha }}-2z_{1}\frac{%
\partial }{\partial a_{\alpha -\delta _{1}}}-z_{2}\frac{\partial }{\partial
a_{\alpha -\delta _{2}}}+2z_{1}z_{2}\frac{\partial }{\partial a_{\alpha
-\delta _{1}-\delta _{2}}} \\
&&+z_{1}^{2}\frac{\partial }{\partial a_{\alpha -2\delta _{1}}}%
-z_{1}^{2}z_{2}\frac{\partial }{\partial a_{\alpha -2\delta _{1}-\delta _{2}}%
}.
\end{eqnarray*}

Similar vector fields are constructed by permuting the z-variables, and
changing the index $\alpha $ as indicated by the permutation. The pole order
of the previous vector fields is equal to 3.

\begin{lemma}
For any $(v_{i})_{1\leq i\leq 3}\in \mathbb{C}^{3},$ there exist $v_{\alpha
}(a),$ with degree at most 1 in the variables $(a_{\gamma }),$ such that $V:=%
\underset{\alpha }{\sum }v_{\alpha }(a)\frac{\partial }{\partial a_{\alpha }}%
+\underset{1\leq j\leq 2}{\sum }v_{j}\frac{\partial }{\partial z_{j}}%
+v_{3}z_{3}\frac{\partial }{\partial z_{3}}$ is tangent to $\overline{%
J_{2}^{v}}(\mathcal{Y}_{0})$ at each point.
\end{lemma}

\begin{proof}
First, we substitute equations 1, 2, 3 in equations 4, 5, 6 to get rid of $%
z_{3},\xi _{3}^{(i)}(1\leq i\leq 2).$ Then, we impose the additional
conditions of vanishing for the coefficients of $\xi _{j}^{(1)}$ in the
second equation (respectively of $\xi _{j}^{(1)}\xi _{k}^{(1)}$ in the third
equation) for any $1\leq j\leq k\leq 2$. Then the coefficients of $\xi
_{j}^{(2)}$ are automatically zero in the third equation. The resulting
equations are
\begin{equation*}
\underset{\left| \alpha \right| \leq d}{\sum }v_{\alpha }z^{\alpha }+%
\underset{j=1}{\overset{2}{\sum }}\underset{\left| \alpha \right| \leq d}{%
\sum }a_{\alpha }\frac{\partial z^{\alpha }}{\partial z_{j}}v_{j}-dv_{3}%
\underset{\left| \alpha \right| \leq d}{\sum }a_{\alpha }z^{\alpha }=0,
\end{equation*}

\begin{equation*}
\underset{\left| \alpha \right| \leq d}{\sum }v_{\alpha }\frac{\partial
z^{\alpha }}{\partial z_{j}}+\underset{k=1}{\overset{2}{\sum }}\underset{%
\left| \alpha \right| \leq d}{\sum }a_{\alpha }\frac{\partial ^{2}z^{\alpha }%
}{\partial z_{j}\partial z_{k}}v_{k}-dv_{3}\underset{\left| \alpha \right|
\leq d}{\sum }a_{\alpha }\frac{\partial z^{\alpha }}{\partial z_{j}}=0,
\end{equation*}

\begin{equation*}
\underset{\left| \alpha \right| \leq d}{\sum }\frac{\partial ^{2}z^{\alpha }%
}{\partial z_{j}\partial z_{k}}v_{\alpha }+\underset{l=1}{\overset{2}{\sum }}%
\underset{\left| \alpha \right| \leq d}{\sum }a_{\alpha }\frac{\partial
^{3}z^{\alpha }}{\partial z_{j}\partial z_{k}\partial z_{l}}v_{l}-dv_{3}%
\underset{\left| \alpha \right| \leq d}{\sum }a_{\alpha }\frac{\partial
^{2}z^{\alpha }}{\partial z_{j}\partial z_{k}}=0.
\end{equation*}

Now we can observe that if the $v_{\alpha }(a)$ satisfy the first equation,
they automatically satisfy the other ones because the $v_{\alpha }$ are
constants with respect to $z$. Therefore it is sufficient to find $%
(v_{\alpha })$ satisfying the first equation. We identify the coefficients
of $z^{\rho }=z_{1}^{\rho _{1}}$ $z_{2}^{\rho _{2}}:$%
\begin{equation*}
v_{\rho }+\underset{j=1}{\overset{2}{\sum }}a_{\rho +\delta _{j}}v_{j}(\rho
_{j}+1)-dv_{3}a_{\rho }=0.
\end{equation*}
\end{proof}

Another family of vector fields can be obtained in the following way.
Consider a $3\times 3$-matrix $A=\left(
\begin{array}{ccc}
A_{1}^{1} & A_{1}^{2} & 0 \\
A_{2}^{1} & A_{2}^{2} & 0 \\
A_{3}^{1} & A_{3}^{2} & 0
\end{array}
\right) \in \mathcal{M}_{3}(\mathbb{C})$ and let $\widetilde{V}:=\underset{%
j,k}{\sum }w_{j}^{(k)}\frac{\partial }{\partial \xi _{j}^{(k)}},$ where $%
w^{(k)}:=A\xi ^{(k)},$ for $k=1,2.$

\begin{lemma}
\label{l1}There exist polynomials $v_{\alpha }(z,a):=\underset{\left| \beta
\right| \leq 2}{\sum }v_{\beta }^{\alpha }(a)z^{\beta }$ where each
coefficient $v_{\beta }^{\alpha }$ has degree at most 1 in the variables $%
(a_{\gamma })$ such that
\begin{equation*}
V:=\underset{\alpha }{\sum }v_{\alpha }(z,a)\frac{\partial }{\partial
a_{\alpha }}+\widetilde{V}
\end{equation*}
is tangent to $\overline{J_{2}^{v}}(\mathcal{Y}_{0})$ at each point.
\end{lemma}

\begin{proof}
First, we substitute equations 1, 2, 3 in equations 4, 5, 6 to get rid of $%
z_{3},\xi _{3}^{(i)}(1\leq i\leq 2).$ We impose the additional conditions of
vanishing for the coefficients of $\xi _{j}^{(1)}$ in the second equation
(respectively of $\xi _{j}^{(1)}\xi _{k}^{(1)}$ in the third equation) for
any $1\leq j\leq k\leq 2$. Then the coefficients of $\xi _{j}^{(2)}$ are
automatically zero in the third equation. The resulting equations are
\begin{equation*}
\underset{\left| \alpha \right| \leq d}{\sum }v_{\alpha }z^{\alpha }=0,\text{
\ }(7)
\end{equation*}

\begin{equation*}
\underset{\left| \alpha \right| \leq d}{\sum }v_{\alpha }\frac{\partial
z^{\alpha }}{\partial z_{j}}+\underset{k=1}{\overset{2}{\sum }}\underset{%
\left| \alpha \right| \leq d}{\sum }a_{\alpha }\frac{\partial z^{\alpha }}{%
\partial z_{k}}A_{k}^{j}-dA_{3}^{j}\underset{\left| \alpha \right| \leq d}{%
\sum }a_{\alpha }z^{\alpha }=0,\text{ \ \ }(8_{j})
\end{equation*}

\begin{equation*}
\underset{\alpha }{\sum }\frac{\partial ^{2}z^{\alpha }}{\partial
z_{j}\partial z_{k}}v_{\alpha }+\underset{\alpha ,p}{\sum }a_{\alpha }\frac{%
\partial ^{2}z^{\alpha }}{\partial z_{j}\partial z_{p}}A_{p}^{k}+\underset{%
\alpha ,p}{\sum }a_{\alpha }\frac{\partial ^{2}z^{\alpha }}{\partial
z_{k}\partial z_{p}}A_{p}^{j}-2dA_{3}^{j}\underset{\alpha }{\sum }a_{\alpha }%
\frac{\partial z^{\alpha }}{\partial z_{k}}=0.\text{ \ }(9_{jk})
\end{equation*}

The equations for the unknowns $v_{\beta }^{\alpha }$ are obtained by
identifying the coefficients of the monomials $z^{\rho }$ in the above
equations.

The monomials $z^{\rho }$ in (7) are $z_{1}^{\rho _{1}}$ $z_{2}^{\rho _{2}}$
with $\sum \rho _{i}\leq d$.

If all the components of $\rho $ are greater than 2, then we obtain the
following system

10. The coefficient of $z^{\rho }$ in (7) impose the condition
\begin{equation*}
\underset{\alpha +\beta =\rho }{\sum }v_{\beta }^{\alpha }=0.
\end{equation*}

11$_{j}.$ The coefficient of the monomial $z^{\rho -\delta _{j}}$ in $%
(8_{j}) $ impose the condition
\begin{equation*}
\underset{\alpha +\beta =\rho }{\sum }\alpha _{j}v_{\beta }^{\alpha
}=l_{j}(a).
\end{equation*}
where $l_{j}$ is a linear expression in the $a$-variables.

12$_{jj}.$ For $j=1,\dots,2$ the coefficient of the monomial $z^{\rho -2\delta
_{j}}$ in $(9_{jj})$ impose the condition
\begin{equation*}
\underset{\alpha +\beta =\rho }{\sum }\alpha _{j}(\alpha _{j}-1)v_{\beta
}^{\alpha }=l_{jj}(a).
\end{equation*}

12$_{jk}.$ For $1\leq j<k\leq 2$ the coefficient of the monomial $z^{\rho
-\delta _{j}-\delta _{k}}$ in $(9_{jk})$ impose the condition
\begin{equation*}
\underset{\alpha +\beta =\rho }{\sum }\alpha _{j}\alpha _{k}v_{\beta
}^{\alpha }=l_{jk}(a).
\end{equation*}

The determinant of the matrix associated to the system is not zero. Indeed,
for each $\rho $ the matrix whose column $C_{\beta }$ consists of the
partial derivatives of order at most 2 of the monomial $z^{\rho -\beta }$
has the same determinant, at the point $z_{0}=(1,1),$ as our system.
Therefore if the determinant is zero, we would have a non-identically zero
polynomial
\begin{equation*}
Q(z)=\underset{\beta }{\sum }a_{\beta }z^{\rho -\beta },
\end{equation*}
such that all its partial derivatives of order less or equal to 2 vanish at $%
z_{0}.$ Thus the same is true for
\begin{equation*}
P(z)=z^{\rho }Q(\frac{1}{z_{1}},\frac{1}{z_{2}})=\underset{\beta }{\sum }%
a_{\beta }z^{\beta }.
\end{equation*}
But this implies $P\equiv 0.$

Finally, we conclude by Cramer's rule. The systems we have to solve are
never over determined. The lemma is proved.
\end{proof}

\begin{remark}
We have chosen the matrix $A$ with this form because we are interested to
prove the global generation statement on $\overline{J_{2}^{v}}(\mathcal{Y}%
)\backslash (\Sigma \cup p^{-1}(H))$ where $\Sigma $ is the closure of $$%
\Sigma _{0}=\{(z,a,\xi ^{(1)},\xi ^{(2)})\in \overline{J_{2}^{v}}(\mathcal{Y}%
_{0}) / \det \left( \xi _{i}^{(j)}\right) _{1\leq i,j\leq 2}=0\}.$$
\end{remark}

\begin{proposition}
The vector bundle $T_{\overline{J_{2}^{v}}(\mathcal{Y})}\otimes \mathcal{O}_{%
\mathbb{P}^{3}}(7)\otimes \mathcal{O}_{\mathbb{P}^{N_{d}+1}}(\ast )$ is
generated by its global sections on $\overline{J_{2}^{v}}(\mathcal{Y}%
)\backslash (\Sigma \cup p^{-1}(H)).$
\end{proposition}

\begin{proof}
From the preceding lemmas, we are reduced to consider $V=\underset{\left|
\alpha \right| \leq 3}{\sum }v_{\alpha }\frac{\partial }{\partial a_{\alpha }%
}.$ The conditions for $V$ to be tangent to $\overline{J_{3}^{v}}(\mathcal{Y}%
_{0})$ are
\begin{equation*}
\underset{\left| \alpha \right| \leq 2}{\sum }v_{\alpha }z^{\alpha }=0,
\end{equation*}

\begin{equation*}
\underset{j=1}{\overset{2}{\sum }}\underset{\left| \alpha \right| \leq 2}{%
\sum }v_{\alpha }\frac{\partial z^{\alpha }}{\partial z_{j}}\xi _{j}^{(1)}=0,
\end{equation*}

\begin{equation*}
\underset{\left| \alpha \right| \leq 2}{\sum }(\underset{j=1}{\overset{2}{%
\sum }}\frac{\partial z^{\alpha }}{\partial z_{j}}\xi _{j}^{(2)}+\underset{%
j,k=1}{\overset{2}{\sum }}\frac{\partial ^{2}z^{\alpha }}{\partial
z_{j}\partial z_{k}}\xi _{j}^{(1)}\xi _{k}^{(1)})v_{\alpha }=0.
\end{equation*}

We have $W_{12}:=\det (\xi _{j}^{(i)})_{1\leq i,j\leq 2}\neq 0.$ Then we can
solve the previous system with $v_{00},v_{10},v_{01}$ as unknowns. By the
Cramer rule, each of the previous quantity is a linear combination of the $%
v_{\alpha },$ $\left| \alpha \right| \leq 2,$ $\alpha \neq (00),$ $(10),$ $%
(01)$ with coefficients rational functions in $z,\xi ^{(1)},\xi ^{(2)}.$ The
denominator is $W_{12}$ and the numerator is a polynomial whose monomials
have either degree at most 2 in $z,$ and at most 1 in $\xi ^{(1)}$ and $\xi
^{(2)},$ or degree 1 in $z$ and three in $\xi ^{(1)}.$

$\xi ^{(1)}$ has a pole of order 2, $\xi ^{(2)}$ has a pole of order 3 therefore the previous vector field has a pole of order at most
7.
\end{proof}

\begin{corollary}
The vector bundle $T_{\overline{J_{2}^{v}}(\mathbb{P}^{2}\times \mathbb{P}%
^{N_{d}})}\otimes \mathcal{O}_{\mathbb{P}^{2}}(7)\otimes \mathcal{O}_{%
\mathbb{P}^{N_{d}}}(\ast )$ is generated by its global sections on $%
\overline{J_{2}^{v}}(\mathbb{P}^{2}\times \mathbb{P}^{N_{d}})\backslash (\pi
_{2}(\Sigma )\cup p_{2}^{-1}(\mathcal{X)})$, where $p_{2}$ is the natural projection $\overline{J_{2}^{v}}(\mathbb{P}^{2}\times \mathbb{P}^{N_{d}}) \rightarrow \mathbb{P}^{2}\times \mathbb{P}^{N_{d}}$.
\end{corollary}

\begin{remark}
The pole order $7$ is the same as in the compact case of \cite{Paun05}.
\end{remark}

\begin{remark}
If the second derivative of $f:(\mathbb{C},0)\rightarrow \mathbb{P}%
^{2}\times \mathbb{P}^{N_{d}}\backslash \mathcal{X}$ lies inside $\pi
_{2}(\Sigma )$ then the image of $f$ is contained in a line. Therefore, as far as we are interested in the algebraic degeneracy of $f$, it is no loss of generality to work away from $\Sigma$.
\end{remark}

\subsection{The two-components case}
By the previous method we obtain the same global generation statement, using the same notations,

\begin{proposition}
The vector bundle $T_{\overline{J_{2}^{v}}(\mathbb{P}^{2}\times \mathbb{P}%
^{N_{d}})}\otimes \mathcal{O}_{\mathbb{P}^{2}}(7)\otimes \mathcal{O}_{%
\mathbb{P}^{N_{d}}}(\ast )$ is generated by its global sections on $%
\overline{J_{2}^{v}}(\mathbb{P}^{2}\times \mathbb{P}^{N_{d}})\backslash (\pi
_{2}(\Sigma )\cup p_{2}^{-1}(\mathcal{X)})$.
\end{proposition}

The proof goes along the same lines considering $$\mathcal{X=X}_{1}\cup \mathcal{X}_{2}\subset \mathbb{P%
}^{2}\times \mathbb{P}^{N_{d_{1}}}\times \mathbb{P}^{N_{d_{2}}},$$ where $%
\mathcal{X}_{i}$ ($1\leq i\leq 2)$ is the universal curve of degree $d_{i}$
given by the equation
\begin{equation*}
\underset{\left| \alpha \right| =d_{i}}{\sum }a_{\alpha }^{(i)}Z^{\alpha }=0,%
\text{ where }[a^{(i)}]\in \mathbb{P}^{N_{d_{i}}}\text{ and }[Z]\in \mathbb{P%
}^{2},
\end{equation*}

and $\mathcal{Y=Y}_{1}\cap \mathcal{Y}_{2}\subset \mathbb{P}%
^{4}\times U$ where
\begin{eqnarray*}
\mathcal{Y}_{1} &=&(a_{000d_{1}0}^{(1)}Z_{3}^{d_{1}}+\underset{\left| \alpha
\right| =d_{1}}{\sum }a_{\alpha }^{(1)}Z^{\alpha }=0)\subset \mathbb{P}%
^{4}\times U, \\
\mathcal{Y}_{2} &=&(a_{0000d_{2}}^{(2)}Z_{4}^{d_{2}}+\underset{\left| \alpha
\right| =d_{2}}{\sum }a_{\alpha }^{(2)}Z^{\alpha }=0)\subset \mathbb{P}%
^{4}\times U
\end{eqnarray*}
where $U$ is the open subset of $\mathbb{P}^{N_{d_{1}}+1}%
\times \mathbb{P}^{N_{d_{2}}+1}$ defined by $$U:=(a_{000d_{1}0}^{(1)}\neq 0)\cap \left( \underset{\left| \alpha
\right| =d_{1},\alpha _{^{\prime }},\alpha _{5}=0}{\cup }(a_{\alpha
}^{(1)}\neq 0)\right) \times (a_{0000d_{2}}^{(2)}\neq 0)\cap \left(
\underset{\left| \alpha \right| =d_{2},\alpha _{^{\prime }},\alpha _{5}=0}{%
\cup }(a_{\alpha }^{(2)}\neq 0)\right).$$ 

Then we apply the previous method to construct meromorphic vector fields on $\overline{J_{2}^{v}}(%
\mathbb{P}^{2}\times \mathbb{P}^{N_{d_{1}}}\times \mathbb{P}^{N_{d_{2}}}).$

\section{Logarithmic jet differentials}

In this section we recall the basic facts about logarithmic jet
differentials following G. Dethloff and S.Lu \cite{DL96}. Let $X$ be a
complex manifold with a normal crossing divisor $D.$

Let $(X,D)$ be the corresponding complex log-manifold. We start with the
directed manifold $(X,\overline{T}_{X})$ where $\overline{T}_{X}=T_{X}(-\log
D).$ We define $X_{1}:=\mathbb{P(}\overline{T}_{X}),$ $D_{1}=\pi ^{\ast }(D)$
and $V_{1}\subset T_{X_{1}}:$%
\begin{equation*}
V_{1,(x,[v])}:=\{\xi \in \overline{T}_{X_{1},(x,[v])}(-\log D_{1})\text{ };%
\text{ }\pi _{\ast }\xi \in \mathbb{C}v\},
\end{equation*}
where $\pi :X_{1}\rightarrow X$ is the natural projection. If $f:(\mathbb{C}%
,0)\rightarrow (X\backslash D,x)$ is a germ of holomorphic curve then it can
be lifted to $X_{1}\backslash D_{1}$ as $f_{[1]}.$

By induction, we obtain a tower of varieties $(X_{k},D_{k},V_{k})$ with $\pi
_{k}:X_{k}\rightarrow X$ as the natural projection. We have a tautological
line bundle $\mathcal{O}_{X_{k}}(1)$ and we denote $u_{k}:=c_{1}(\mathcal{O}%
_{X_{k}}(1)).$

Let us consider the direct image $\pi _{k\ast }(\mathcal{O}_{X_{k}}(m)).$
It is a locally free sheaf denoted $E_{k,m}\overline{T}_{X}^{\ast }$
generated by all polynomial operators in the derivatives of order $1,2,...,k$
of $f$, together with the extra function $\log s_{j}(f)$ along the $j-th$
component of $D,$ which are moreover invariant under arbitrary changes of
parametrization: a germ of operator $Q\in E_{k,m}\overline{T}_{X}^{\ast }$
is characterized by the condition that, for every germ  of holomorphic curve $f:(\mathbb{C}%
,0)\rightarrow (X\backslash D,x)$
and every germ $\phi \in $ $\mathbb{G}_{k}$ of $k$-jet biholomorphisms of $(%
\mathbb{C},0),$%
\begin{equation*}
Q(f\circ \phi )=\phi ^{\prime m}Q(f)\circ \phi .
\end{equation*}

\bigskip

The following theorem makes clear the use of jet differentials in the study
of hyperbolicity:

\bigskip

\begin{theorem}\label{tvan}(\cite{GG80}, \cite{De95}, \cite{DL96}). Assume
that there exist integers $k,m>0$ and an ample line bundle $L$ on X such that
\begin{equation*}
H^{0}(X_{k},\mathcal{O}_{X_{k}}(m)\otimes \pi _{k}^{\ast }L^{-1})\simeq
H^{0}(X,E_{k,m}\overline{T}_{X}^{\ast }\otimes L^{-1})
\end{equation*}
has non zero sections $\sigma _{1},\dots,\sigma _{N}$. Let $Z\subset X_{k}$ be the base locus of these sections. Then every
entire curve $f:\mathbb{C}\rightarrow X\backslash D$ is such that
$f_{[k]}(\mathbb{C})\subset Z$. In other words, for every global $%
\mathbb{G}_{k}-$ invariant polynomial differential operator P with
values in $L^{-1}$, every entire curve $f:\mathbb{C}\rightarrow
X\backslash D$ must satisfy the algebraic differential equation $%
P(f)=0.$
\end{theorem}

\bigskip

In the case of logarithmic surfaces $(X,D)$, we have the following
filtration (see \cite{E.G}) of log-jet differentials of order 2:
\begin{equation*}
Gr^{\bullet }E_{2,m}\overline{T}_{X}^{\ast }=\underset{0\leq j\leq m/3}{%
\oplus }S^{m-3j}\overline{T}_{X}^{\ast }\otimes \overline{K}_{X}^{\otimes j}.
\end{equation*}

A Riemann-Roch calculation based on the above filtration yields
\begin{equation*}
\chi (X,E_{2,m}\overline{T}_{X}^{\ast })=\frac{m^{4}}{648}(13\overline{c}%
_{1}^{2}-9\overline{c}_{2})+O(m^{3}),
\end{equation*}
where $\overline{c}_{1}$ and $\overline{c}_{2}$ denote the logarithmic Chern classes. This gives by Bogomolov's vanishing theorem \cite{Bo79}:

\begin{theorem}
\cite{E.G} If $(X,D)$ is an algebraic log surface of log general type and A
an ample line bundle over X, then
\begin{equation*}
h^{0}(X,E_{2,m}\overline{T}_{X}^{\ast }\otimes \mathcal{O}(-A))\geq \frac{%
m^{4}}{648}(13\overline{c}_{1}^{2}-9\overline{c}_{2})+O(m^{3}).
\end{equation*}
\end{theorem}

\begin{corollary}
\cite{E.G} Let $C\subset \mathbb{P}^{2}$ be a smooth curve of degree $d\geq
11$ and A
an ample line bundle. Then $h^{0}(X,E_{2,m}\overline{T}_{X}^{\ast }\otimes \mathcal{O}%
(-A))\neq 0$ for $m$ large enough.
\end{corollary}

\begin{corollary}
\cite{Rou03} Let $C=C_{1}\cup C_{2}$ be a normal crossing complex algebraic curve in $\mathbb{P}%
^{2}$ having two irreducible smooth components $C_{1}$ and $C_{2}$ of degrees $d_{1}\leq d_{2}$ and A
an ample line bundle. Then $h^{0}(X,E_{2,m}\overline{T}_{X}^{\ast }\otimes \mathcal{O}%
(-A))\neq 0$ for $m$ large enough if the degrees satisfy
$$\begin{array}{cccc}
either & d_{1} \geq 3, \\
or & d_{1} =2 & and & d_{2}\geq 5,\\
or & d_{1}=1 &and& d_{2}\geq 7.
\end{array}
$$
\end{corollary}

\bigskip

\section{Proof of theorem \ref{t1}}

Let us consider an entire curve $f:\mathbb{C}\rightarrow \mathbb{P}%
^{2}\backslash C$ for a generic curve in $\mathbb{P}^{2}$ of degree $d\geq
14.$ Let us assume that the projectivized first derivative $f_{[1]}:\mathbb{%
C} \rightarrow X_{1}$ is Zariski dense. By the proposition of the previous
section we have a section
\begin{equation*}
\sigma \in H^{0}(\mathbb{P}^{2},E_{2,m}\overline{T}_{\mathbb{P}^{2}}^{\ast
}\otimes \overline{K}_{\mathbb{P}^{2}}^{-t})\simeq H^{0}((\mathbb{P}%
^{2})_{2},\mathcal{O}_{(\mathbb{P}^{2})_{2}}(m)\otimes \pi _{2}^{\ast }%
\overline{K}_{\mathbb{P}^{2}}^{-t}).
\end{equation*}
with zero set $Z$ and vanishing order $t(d-3).$ Consider the family
\begin{equation*}
\mathcal{X}\subset \mathbb{P}^{2}\times \mathbb{P}^{N_{d}}
\end{equation*}
of curves of degree $d$ in $\mathbb{P}^{2}.$ General semicontinuity
arguments concerning the cohomology groups show the existence of a Zariski
open set $U_{d}\subset \mathbb{P}^{N_{d}}$ such that for any $a\in U_{d},$
there exists an irreducible and reduced divisor
\begin{equation*}
Z_{a}=(P_{a}=0)\subset (\mathbb{P}_{a}^{2})_{2},
\end{equation*}
where
\begin{equation*}
P_{a}\in H^{0}((\mathbb{P}_{a}^{2})_{2},\mathcal{O}_{(\mathbb{P}%
_{a}^{2})_{2}}(m)\otimes \pi _{2}^{\ast }\overline{K}_{(\mathbb{P}%
_{a}^{2})}^{-t}),
\end{equation*}
such that the family $(P_{a})_{a\in U_{d}}$ varies holomorphically. We have
the following numerical criterion due to El Goul

\begin{proposition}
\cite{E.G} Let $(X,D)$ be a log surface of log general type with $Pic(X)=%
\mathbb{Z}.$ Suppose that
\begin{equation*}
m(13\overline{c}_{1}^{2}-9\overline{c}_{2})>12t\overline{c}_{1}^{2},
\end{equation*}
then there exists a divisor $Y_{1}\subset X_{1}$ such that $%
im(f_{[1]})\subset Y_{1}.$
\end{proposition}

Therefore since we assume that the projectivized first derivative $f_{[1]}:%
\mathbb{C} \rightarrow X_{1}$ is Zariski dense we obtain the following
estimate for the vanishing order
\begin{equation*}
t\geq \frac{m(13\overline{c}_{1}^{2}-9\overline{c}_{2})}{12\overline{c}%
_{1}^{2}}.
\end{equation*}

Now we consider $P$ as a holomorphic function on $\overline{J_{2}^{v}}(%
\mathbb{P}^{2}\times \mathbb{P}^{N_{d}})_{U_{d}}$ and differentiate it with
the meromorphic vector fields constructed before. Take $v\in H^{0}(\overline{%
J_{2}^{v}}(\mathbb{P}^{2}\times \mathbb{P}^{N_{d}})_{U_{d}},T_{\overline{%
J_{2}^{v}}(\mathbb{P}^{2}\times \mathbb{P}^{N_{d}})_{U_{d}}}\otimes \mathcal{%
O}_{\mathbb{P}^{2}}(3))$ such a vector field, then the restriction of $dP(v)$
to $Z_{a}$ is a section of the bundle $$\mathcal{O}_{(\mathbb{P}%
_{a}^{2})_{2}}(m)\otimes \mathcal{O}_{\mathbb{P}_{a}^{2}}(3-t(d-3))\left|
Z_{a}\right. .$$

From the previous proposition, we have
\begin{equation*}
t(d-3)\geq \frac{m(13\overline{c}_{1}^{2}-9\overline{c}_{2})}{12\overline{c}%
_{1}^{2}}(d-3)>3,
\end{equation*}
if $d>14.$

Therefore if $f_{[1]}:\mathbb{C} \rightarrow X_{1}$ is Zariski dense, we have
that the restriction of $dP(v)$ to $Z_{a}$ must vanish. Then there exists a
section $$\rho _{v}\in H^{0}(\mathbb{P}_{a}^{2},\mathcal{O}_{\mathbb{P}%
_{a}^{2}}(3+t(d-3)),$$ such that
\begin{equation*}
dP(v)=\rho _{v}P.
\end{equation*}
As in the compact case, we have the following proposition:

\begin{proposition}
Let $f:\mathbb{C}\rightarrow \mathbb{P}^{2}\backslash C$ such that $f_{[1]}:%
\mathbb{C} \rightarrow X_{1}$ is Zariski dense. Then the weighted degree $m$
of the algebraic family of sections $(P_{a})$ verifies $m\geq 6$.
\end{proposition}

\begin{proof}
Following \cite{E.G}, we have for $m\leq 5$ the exact sequence
\begin{equation*}
0\rightarrow S^{m}\overline{T}_{\mathbb{P}^{2}}^{\ast }\otimes \overline{K}_{%
\mathbb{P}^{2}}^{-t}\rightarrow E_{2,m}\overline{T}_{\mathbb{P}^{2}}^{\ast
}\otimes \overline{K}_{\mathbb{P}^{2}}^{-t}\overset{\phi }{\rightarrow }%
S^{m-3}\overline{T}_{\mathbb{P}^{2}}^{\ast }\otimes \overline{K}_{\mathbb{P}%
^{2}}^{1-t}\rightarrow 0,
\end{equation*}
which gives an injective morphism for any positive $t$%
\begin{equation*}
\Phi :H^{0}(\mathbb{P}^{2},E_{2,m}\overline{T}_{\mathbb{P}^{2}}^{\ast
}\otimes \overline{K}_{\mathbb{P}^{2}}^{-t})\rightarrow H^{0}(\mathbb{P}%
^{2},S^{m-3}\overline{T}_{\mathbb{P}^{2}}^{\ast }\otimes \overline{K}_{%
\mathbb{P}^{2}}^{1-t}),
\end{equation*}
because $H^{0}(\mathbb{P}^{2},S^{m}\overline{T}_{\mathbb{P}^{2}}^{\ast
}\otimes \overline{K}_{\mathbb{P}^{2}}^{-t})=0$ (see \cite{E.G}).

Let us consider as above the logarithmic manifold $(\mathcal{Y},H)$ with the
log-morphism $\pi :(\mathcal{Y},H)\rightarrow (\mathbb{P}^{2}\times \mathbb{P%
}^{N_{d}},\mathcal{X}).$ We can take the pull-back with $\pi $ of the
algebraic family of sections $(P_{a})$ providing logarithmic 2-jet
differentials on $\mathcal{Y}_{a}.$ Let us do some local computations on
these 2-jet differentials. We take some affine coordinates on $\mathbb{C}%
^{3} $ and the equation of $\mathcal{Y}_{0}$ is
\begin{equation*}
\mathcal{Y}_{0}:=(z_{3}^{d}+\underset{\alpha }{\sum }a_{\alpha }z^{\alpha
}=0)\subset \mathbb{C}^{3}\times U,
\end{equation*}
with $H_{0}:=(z_{3}=0).$

Now, the proof is similar to the compact case (see \cite{Paun05}) which we
recall for the convenience of the reader.

If we assume that $\underset{\alpha }{\sum }a_{\alpha }\frac{\partial
z^{\alpha }}{\partial z_{1}}\neq 0$ identically on $\mathcal{Y}_{0},$ we can
write a logarithmic 2-jet differential on the corresponding affine open set
with logarithmic coordinates in the following way
\begin{equation*}
Q(z,a,\xi ^{(1)},\xi ^{(2)})=R_{0}(z,a,\xi ^{(1)})(\xi _{2}^{(1)}\xi
_{3}^{(2)}-\xi _{3}^{(1)}\xi _{2}^{(2)})+R_{1}(z,a,\xi ^{(1)}),
\end{equation*}
where $R_{0}$ and $R_{1}$ are local symmetric differentials, of degree $m-3$
and $m.$ For a generic point $z\in \mathcal{Y}_{0},$ $R_{0}$ and $R_{1}$ are
not identically zero since the zero set of the $(P_{a})$ is irreducible.

For generic $z_{0}=(z_{1}^{0},z_{2}^{0},z_{3}^{0})\in \mathcal{Y}_{0},$ $%
R_{0}$ and $R_{1}$ are not identically zero. By translation we can assume
that $z_{1}^{0}=z_{2}^{0}=0.$ We can make a linear transformation on $%
(z_{1},z_{2})$ to diagonalize the quadratic part and the equation of $%
\mathcal{Y}_{0}$ becomes
\begin{equation*}
\mathcal{Y}_{0}=(z_{3}^{d}+\underset{3\leq \alpha \leq d}{\sum }c_{\alpha
}z^{\alpha }+c_{0}(z_{1}^{2}+z_{2}^{2}+c_{100}z_{1})=0).
\end{equation*}

Notice that the equation of the divisor $H_{0}$ is still $(z_{3}=0).$ In
these coordinates, we consider as above the manifold $J_{2}^{v}(\mathcal{Y}%
_{0})$ and $\widetilde{V}:=\underset{j}{\sum }\xi _{2}^{(j)}\frac{\partial }{%
\partial \xi _{3}^{(j)}}-\xi _{3}^{(j)}\frac{\partial }{\partial \xi
_{2}^{(j)}}.$ According to lemma \ref{l1}, there exists a global meromorphic
vector field $V\in H^{0}(J_{2}^{v}(\mathcal{Y}),T_{J_{2}^{v}(\mathcal{Y}%
)}\otimes \mathcal{O}_{\mathbb{P}^{3}}(3))$ such that $V:=\underset{\alpha }{%
\sum }v_{\alpha }\frac{\partial }{\partial a_{\alpha }}+\widetilde{V}.$

Moreover, from the proof of the lemma we see that for each index $\alpha
,\left| \alpha \right| \leq 2,$ we have $v_{\alpha }(0,0)=0.$

Now, using the first package of vector fields and the relation $dP(v)=\rho
_{v}P$, we obtain that there exists $\lambda \in \mathbb{C}$ such that
\begin{equation*}
\left( \underset{j}{\sum }\xi _{2}^{(j)}\frac{\partial }{\partial \xi
_{3}^{(j)}}-\xi _{3}^{(j)}\frac{\partial }{\partial \xi _{2}^{(j)}}\right)
Q(z_{0},a_{0},\xi ^{(1)},\xi ^{(2)})=\lambda Q(z_{0},a_{0},\xi ^{(1)},\xi
^{(2)}),
\end{equation*}
for any 2-jet $(\xi ^{(1)},\xi ^{(2)})$ of $\mathcal{Y}_{0}$ at $z_{0}.$ We
remark that $$\left( \underset{j}{\sum }\xi _{2}^{(j)}\frac{\partial }{%
\partial \xi _{3}^{(j)}}-\xi _{3}^{(j)}\frac{\partial }{\partial \xi
_{2}^{(j)}}\right) (\xi _{2}^{(1)}\xi _{3}^{(2)}-\xi _{3}^{(1)}\xi
_{2}^{(2)})=0.$$

If $m=3,$ then we have $\lambda =0$ since $R_{0}$ is not $0,$ so
\begin{equation*}
\left( \xi _{2}^{(1)}\frac{\partial }{\partial \xi _{3}^{(1)}}-\xi _{3}^{(1)}%
\frac{\partial }{\partial \xi _{2}^{(1)}}\right) \left( \underset{j=0}{%
\overset{3}{\sum }}R_{j}^{1}\left( \xi _{2}^{(1)}\right) ^{j}\left( \xi
_{3}^{(1)}\right) ^{3-j}\right) =0,
\end{equation*}
which implies that $R_{j}^{1}=0$ for all $j.$

If $m=4,$ then
\begin{equation*}
\left( \xi _{2}^{(1)}\frac{\partial }{\partial \xi _{3}^{(1)}}-\xi _{3}^{(1)}%
\frac{\partial }{\partial \xi _{2}^{(1)}}\right) \left( R_{1}^{0}\xi
_{2}^{(1)}+R_{2}^{0}\xi _{3}^{(1)}\right) =\lambda \left( R_{1}^{0}\xi
_{2}^{(1)}+R_{2}^{0}\xi _{3}^{(1)}\right),
\end{equation*}
and
\begin{equation*}
\left( \xi _{2}^{(1)}\frac{\partial }{\partial \xi _{3}^{(1)}}-\xi _{3}^{(1)}%
\frac{\partial }{\partial \xi _{2}^{(1)}}\right) \left( \underset{j=0}{%
\overset{4}{\sum }}R_{j}^{1}\left( \xi _{2}^{(1)}\right) ^{j}\left( \xi
_{3}^{(1)}\right) ^{3-j}\right) =\lambda \left( \underset{j=0}{\overset{4}{%
\sum }}R_{j}^{1}\left( \xi _{2}^{(1)}\right) ^{j}\left( \xi
_{3}^{(1)}\right) ^{3-j}\right),
\end{equation*}
which imposes the two equations
\begin{eqnarray*}
\lambda ^{2}-1 &=&0, \\
\lambda (\lambda ^{4}+8\lambda +16) &=&0.
\end{eqnarray*}
They do not have common solutions.

If $m=5$ we obtain
\begin{eqnarray*}
\lambda (\lambda ^{2}+4) &=&0, \\
\lambda ^{2}(13+\lambda ^{2})^{2}+9(5+\lambda ^{2}) &=&0,
\end{eqnarray*}
with no common solutions.
\end{proof}

So, with $f:\mathbb{C}\rightarrow \mathbb{P}^{2}\backslash C$ such that $%
f_{[1]}:\mathbb{C} \rightarrow X_{1}$ is Zariski dense we have degree at
least 6 for the jet differential. But, we have proved that $T_{\overline{%
J_{2}^{v}}(\mathbb{P}^{2}\times \mathbb{P}^{N_{d}})}\otimes \mathcal{O}_{%
\mathbb{P}^{2}}(7)_{\left| X_{2}\right. }$ is globally generated on $%
X_{2}\backslash (\pi _{2}(\Sigma )\cup \pi _{2,0}^{\ast }(C)).$ So, we can
find a vector field $v$ such that $dP(v)$ is a \ holomorphic jet
differential vanishing on ample divisor and algebraically independent of $P$
provided that
\begin{equation*}
\frac{m(13\overline{c}_{1}^{2}-9\overline{c}_{2})}{12\overline{c}_{1}^{2}}%
(d-3)>7,
\end{equation*}
which is the case for $d\geq 14$ as $m\geq 6.$ Then we obtain a contradiction to the fact that  $%
f_{[1]}:\mathbb{C} \rightarrow X_{1}$ is Zariski dense. Indeed, if $f_{[2]}$ does not lie in the singular set of $(P=0)$ then for some $t \in \mathbb{C}$, $dP_{f_{[2]}(t)}(v)\neq 0$, which contradicts theorem \ref{tvan}. Therefore $f_{[2]}$ lies in this singular set and $f_{[1]}$ is algebraically degenerate. Then we use El Goul's generalization of McQuillan's results on foliations of surfaces to the logarithmic setting in \cite{E.G} which tells us that the entire curve $f:%
\mathbb{C}\rightarrow \mathbb{P}^{2}\backslash C$ itself is algebraically
degenerate. Finally, algebraic hyperbolicity of
$\mathbb{P}^{2}\backslash C$ (see \cite{PaRou} or \cite{ch01}) implies that $f$ is constant and $\mathbb{P}^{2}\backslash C$ is hyperbolic and hyperbolically embedded in $ \mathbb{P}^{2}$.

\section{Proof of theorem \ref{t2}}

The proof goes along the same lines showing the algebraic degeneracy of the
curve $f:\mathbb{C}\rightarrow \mathbb{P}^{2}\backslash C_{1}\cup C_{2}$
provided that the same numerical condition
\begin{equation*}
\frac{m(13\overline{c}_{1}^{2}-9\overline{c}_{2})}{12\overline{c}_{1}^{2}}%
(d-3)>7,
\end{equation*}
is satisfied.
An easy computation shows that it is the case if either $d_{1}\geq4,$ or $d_{1}=3$ and $d_{2}\geq 5$, or $d_{1}=2$ and $d_{2}\geq
8$, or $d_{1}=1$ and $d_{2}\geq 11$. The hyperbolicity is deduced
from the algebraic hyperbolicity of $\mathbb{P}^{2}\backslash C_{1}\cup
C_{2} $ (see \cite{ch01}).

\noindent \texttt{rousseau@math.u-strasbg.fr}

\noindent D\'{e}partement de Math\'{e}matiques,

\noindent IRMA,\newline
Universit\'{e} Louis Pasteur,

\noindent 7, rue Ren\'{e} Descartes,\newline
\noindent 67084 STRASBOURG CEDEX

\noindent FRANCE


\bigskip

\begin{thebibliography}{99}
\bibitem{BeDu} Berteloot F., Duval J., \textit{Sur l'hyperbolicit\'e de certains compl\'ementaires},  Enseign. Math. II. Ser. 47 (2001), p. 253--267.

\bibitem{Bo79}  Bogomolov F.A., \textit{Holomorphic tensors and vector
bundles on projective varieties}, Math. USSR Izvestija 13 (1979), 499--555.

\bibitem{ch01}  Chen X. , \textit{On Algebraic Hyperbolicity of Log Varieties%
}, Commun. Contemp. Math. 6 (2004), no. 4, 513--559. Also available
as preprint math.AG/0111051.

\bibitem{Cle}  Clemens H., \textit{Curves on generic hypersurface}, Ann.
Sci. Ec. Norm. Sup., 19 (1986), 629--636.

\bibitem{De95}  Demailly J.-P., \textit{Algebraic criteria for Kobayashi
hyperbolic projective varieties and jet differentials}, Proc. Sympos. Pure
Math., vol.62, Amer. Math.Soc., Providence, RI (1997), 285--360.

\bibitem{DEG00}  Demailly J.-P., El Goul J., \textit{Hyperbolicity of
generic surfaces of high degree in projective 3-space}, Amer. J. Math 122 (2000), 515--546.

\bibitem{DL96}  Dethloff G., Lu. S, \textit{Logarithmic jet bundles and
applications}, Osaka J. Math. 38 (2001), 185--237.

\bibitem{DL04}  Dethloff G., Lu. S, \textit{Logarithmic surfaces and
hyperbolicity}, Ann. Inst. Fourier (Grenoble)  57  (2007),  no. 5, 1575--1610.

\bibitem{DSW1}  Dethloff G., Schumacher G., Wong P.M., \textit{Hyperbolicity of the complements of plane algebraic curves},  Amer. J. Math. 117 (1995), p. 573--599.

\bibitem{DSW2}  Dethloff G., Schumacher G., Wong P.M., \textit{On the
hyperbolicity of the complements of curves in algebraic surfaces: the three
component case}, Duke Math. J., vol. 78 (1995), 193--212.

\bibitem{Ein}  Ein L., \textit{Subvarieties of generic complete intersections%
}, Invent. Math. 94 (1988), 163--169.

\bibitem{E.G}  El Goul J., \textit{Logarithmic Jets and Hyperbolicity},
Osaka J. Math. 40 (2003), 469--491.

\bibitem{GG80}  Green M., Griffiths P., \textit{Two applications of
algebraic geometry to entire holomorphic mappings}, The Chern Symposium
1979, Proc. Inter. Sympos. Berkeley, CA, 1979, Springer-Verlag, New-York (1980), 41--74.

\bibitem{Ko70}  Kobayashi S., \textit{Hyperbolic manifolds and holomorphic
mappings}, Marcel Dekker, New York, 1970.

\bibitem{Ko98}  Kobayashi S., \textit{Hyperbolic complex spaces}, Grundlehren der Mathematischen Wissenschaften, 318, Springer-Verlag, Berlin, 1998.

\bibitem{No86}  Noguchi J., \textit{Logarithmic jet spaces and extensions of
de Franchis' Theorem}, Contributions to Several Complex Variables
(Conference in Honor of W. Stoll, Notre Dame 1984), Aspects of
Math., Vieweg, Braunschweig, (1986), 227--249.

\bibitem{NWY} Noguchi J., Winkelmann J., Yamanoi K. \textit{Degeneracy of holomorphic curves into algebraic varieties}, J. Math. Pures Appl. (9)  88  (2007),  no. 3, 293--306.

\bibitem{PaRou}  Pacienza G., Rousseau E., \textit{On the logarithmic
Kobayashi conjecture},  J. Reine Angew. Math.  611  (2007), 221--235.

\bibitem{Paun05}  Paun M., \textit{Vector fields on the total space of
hypersurfaces in the projective space and hyperbolicity},  Math. Ann.  340  (2008),  no. 4, 875--892.

\bibitem{Rou03}  Rousseau E., \textit{Hyperbolicit\'{e} du
compl\'{e}mentaire d'une courbe : le cas de deux composantes}, CRAS Ser. I
336 (2003), 635--640.

\bibitem{Rou05}  Rousseau E., \textit{Etude des jets de Demailly-Semple en
dimension 3}, Ann. Inst. Fourier 56 (2006), 397--421.

\bibitem{Rou2}  Rousseau E., \textit{Equations diff\'{e}rentielles sur les
hypersurfaces de $\mathbb{P}^{4}$, }J. Math. Pures Appl. 86 (2006), 322--341.

\bibitem{Rou06}  Rousseau E., \textit{Weak analytic hyperbolicity of generic
hypersurfaces of high degree in $\mathbb{P}^{4}$},  Ann. Fac. Sci. Toulouse Math. (6)  16  (2007),  no. 2, 369--383.

\bibitem{Rou06bis}  Rousseau E., \textit{Weak analytic hyperbolicity of
complements of generic surfaces of high degree in projective 3-space},
 Osaka J. Math.  44  (2007),  no. 4, 955--971.

\bibitem{SY04}  Siu Y.-T., \textit{Hyperbolicity in complex geometry}, The
legacy of Niels Henrik Abel, Springer, Berlin (2004), 543--566.

\bibitem{Voi}  Voisin C., \textit{On a conjecture of Clemens on rational
curves on hypersurfaces}, J. Diff. Geom., 44 (1996), 200--213.
\end{thebibliography}
\end{document}